\documentclass[10pt]{amsart}

\title{Admissible Dirichlet Series}

\author{Stanley Burris}
\thanks{Research of the first author was supported by a grant from NSERC}
\address{Dept.\ of Pure Mathematics, University of Waterloo,
Waterloo, Ont., Canada N2L 3G1}
\email{snburris@thoralf.uwaterloo.ca}
\author{Karen Yeats}
\address{Dept.\ of Pure Mathematics, University of Waterloo,
Waterloo, Ont., Canada N2L 3G1}
\email{kayeats@uwaterloo.ca}
\thanks{Research of the second author was supported by a NSERC Undergraduate 
Student Research Award}
\thanks{This paper benefited from discussions with R.~Warlimont
after he had received our preliminary version \cite{prog}.
He suggested the use of an integral condition in the definition of 
admissible---the definition continued to evolve after this suggestion. 
Also he found that $\exp\big(\zeta(s)^k\big)$ are examples of admissible 
functions.}

\date{\today}


\newtheorem{theorem}{Theorem}
\newtheorem{corollary}[theorem]{Corollary}
\newtheorem{lemma}[theorem]{Lemma}

\theoremstyle{definition}
\newtheorem{definition}[theorem]{Definition}
\newtheorem{example}[theorem]{Example}
\newtheorem{problem}{Problem}

\newtheorem{remark}[theorem]{Remark}

\newcommand{\Fbox}[1]{\mbox{#1}}

\newcommand{\fa}{\mathbf{a}}
\newcommand{\fb}{\mathbf{b}}
\newcommand{\fc}{\mathbf{c}}

\newcommand{\RV}{\mathsf{RV}}

\newcommand{\bF}{\mathbf{F}}

\newcommand{\bG}{\mathbf{G}}
\newcommand{\bH}{\mathbf{H}}

\newcommand{\sR}{\mathsf{R}}
\newcommand{\bbR}{\mathbb{R}}
\newcommand{\bS}{\mathbf{S}}

\newcommand{\BF}{{\widehat{F}}}
\newcommand{\BG}{{\widehat{G}}}

\newcommand{\ro}{{\rm o}}
\newcommand{\rO}{{\rm O}}

\newcommand{\baralpha} {\quad\text{ as } \sigma \rightarrow \mbox{$\alpha\scriptsize +$}}
\newcommand{\ralpha} {\sigma \rightarrow\alpha\scriptsize +}
\newcommand{\bralpha} {\sigma \rightarrow \mbox{$\alpha\scriptsize +$}}
\newcommand{\domx}{ > 0}

\newcounter{thlistctr}
\newenvironment{thlist}{\ 
\begin{list}%
{\alph{thlistctr}}%
{\setlength{\labelwidth}{2ex}%
\setlength{\labelsep}{1ex}%
\setlength{\leftmargin}{6ex}%
\usecounter{thlistctr}}}%
{\end{list}}

\begin{document}

\begin{abstract}
We propose a definition of admissible Dirichlet series
as the analog of Hayman's \cite{h} 1956 definition of admissible power series.
\end{abstract}

\maketitle

\section{Introduction}

In 1956 Hayman defined  {\em admissible functions}---they are
analytic in a neighborhood of 0 and one can use the saddle point
method to estimate the coefficients of the power series expansion of
such functions. They include the functions $e^z$ and
$\exp\dfrac{1}{1-z}$, are closed under product (of series with the same
radius of convergence) and under exponentiation.

In this paper a notion of admissibility for functions that have Dirichlet series 
expansions is proposed.  We believe that this is a viable analog of Hayman's 
definition because 
(1) this notion of admissible generalizes the conditions of Tenenbaum 
in \cite{te}, 
(2) there is a fundamental theorem (Theorem \ref{main theorem}) that is the analog
of Hayman's fundamental theorem, and
(3) a product of admissible Dirichlet series (with the same abscissa of convergence) 
is again admissible.

\section{Definition of Admissible}

\begin{theorem}\label{apos}
 Suppose the function $\bF(s)$ 
 \begin{thlist}
 \item[\bf A1] 
 has a Dirchlet series expansion
 $\bF(s) = \sum_{n \geq 1}f(n)n^{-s}$, 
 where the coefficients $f(n)$ are nonnegative real, with $f(1)>0$, 
 \item[\bf A2]
 has abscissa of (absolute) convergence $\alpha \in[0,\infty)$, and 
 \item[\bf A3]
 $\bF(s)$ has no zeros in its halfplane of convergence.
 \end{thlist}
 Then there exists a Dirichlet series $\bH(s)$ with real coefficients
 such that 
 $\bF(s) = e^{\bH(s)}$ for $\sigma > \alpha$ 
 where $s = \sigma + it$.
 \end{theorem}
 \begin{proof}
 This is a slight specialization of Theorem 11.14 in Apostol
 \cite{a}.
 \end{proof}

Let $\bbR$ be the set of real numbers.
Assuming that $\bF(s)$ satisfies (A1)--(A3)
and $\bF(s)=e^{\bH(s)}$,
we will need the basic facts about a Taylor series expansion
with remainder of $\bH(\sigma+it)$ about $t=0$.
With
\begin{equation}\label{H derivs}
 \fa(s) \, :=\,  \bH'(s)
 \quad\fb(s) \,:=\, \bH''(s)
 \quad\fc(s) \,:=\, \bH'''(s)
\end{equation}
we have, for $\sigma>\alpha$ and $t\in \bbR$,
\begin{eqnarray}
\bH(\sigma+it)
& =& \bH(\sigma)\,+\,i\fa(\sigma)t\,-\,\frac{\fb(\sigma)}{2}t^2 \,+\, 
\sR(\sigma+it)\label{H rem}
\end{eqnarray}
with the remainder term given by
\begin{eqnarray}
\sR(\sigma+it) 
&=&
-\frac{i}{2}\int_0^t \fc(\sigma+iv)(t - v)^2dv.
\label{rem term}
\end{eqnarray}
Thus we can write
\begin{eqnarray}
\frac{\bF\big(\sigma+it\big)}{\bF(\sigma)}
&=&\exp\Big(i\fa(\sigma)t -\frac{\fb(\sigma)}{2}t^2 + 
\sR\big(\sigma + it\big)\Big)\label{quotient}\\
\frac{\Big|\bF\big(\sigma + it\big)\Big|}{\bF(\sigma)}
&=&\exp\Big(-\frac{\fb(\sigma)}{2}t^2\Big) \cdot
\Big|\exp\Big(\sR\big(\sigma +it\big)\Big)\Big|.\label{abs quotient}
\end{eqnarray}

 \begin{definition}\label{def adm}
 Suppose $\bF(s)$ satisfies (A1)--(A3), with $\bF(s) = e^{\bH(s)}$. 
Let $\fa(s)$ and $\fb(s)$ be the first two
derivatives of $\bH(s)$ as in \eqref{H derivs},
let $\sR(\sigma+it)$ be the remainder term of the Taylor expansion 
for $\bH(\sigma+it)$ as
in \eqref{H rem}, and suppose
  there is a function
  $\delta: (\alpha, \beta) \rightarrow (0, 1)$, for some $\beta>\alpha$, 
such that as $\bralpha$
 \begin{thlist}
 \item[\bf A4]
  $\delta(\sigma) \ \rightarrow\ 0$
 \item[\bf A5]
  $\sigma^2\fb(\sigma) \ \rightarrow\ \infty$
 \item[\bf A6]
$\displaystyle 
\fb(\sigma)\cdot\exp\Big(-\fb(\sigma)\delta(\sigma)^2\Big)\ \rightarrow\ 0$
 \item[\bf A7] 
 $\sR(\sigma+it)\ \rightarrow\ 0
\qquad\text{uniformly for } |t|\le\delta(\sigma)$
 \item[\bf A8]
 $\displaystyle
 \frac{\sigma\sqrt{\fb(\sigma)}}{\bF(\sigma)}
 \int_{|t|\ge\delta(\sigma)} \big|\bF(\sigma+it)\big|\;
\frac{dt}{\sigma^2+ t^2}\ \rightarrow \ 0\;$.
 \end{thlist}
 Then we say $\bF(s)$ is {\bf admissible}, as {\bf witnessed} by $\delta(\sigma)$.
 \end{definition}

\begin{remark} Except for $\S$\ref{closure} we can replace (A6) and (A8)
by the following, giving a more general notion of admissible: 
 \begin{thlist}
 \item[\bf A6-]
  $\displaystyle \fb(\sigma)\delta(\sigma)^2\ \rightarrow\ \infty$
 \item[\bf A8-]
 $\displaystyle
 \frac{\sigma\sqrt{\fb(\sigma)}}{\bF(\sigma)}
 \int_{|t|\ge\delta(\sigma)} \bF(\sigma+it)x^{it}
\frac{dt}{(\sigma+it)(\sigma+1+it)}\ \rightarrow \ 0
\\\text{uniformly for } x>0$
 \end{thlist}
The full strength of (A6) and (A8) are used to prove the
product theorem in $\S$\ref{closure}.
\end{remark}

\section{Asymptotic Estimates and Regular Variation}

In this section we assume that $\bF(s)$ is admissible,
witnessed by $\delta(\sigma)$. 
(A5) implies
\begin{equation}\label{big b}
\fb(\sigma)\ \rightarrow\ \infty \qquad \text{as } \sigma\rightarrow
\mbox{$\alpha+$},
\end{equation}
so there is a $\beta>\alpha$ such that
\begin{equation}
\fb(\sigma)>0\qquad \text{for } \sigma\in(\alpha,\beta).
\end{equation}
We will consistently use $\beta$ as a number in $(\alpha,\infty)$
such that $\fb(\sigma)>0$ on $(\alpha,\beta)$, keeping the original requirement
that $\delta(\sigma)$ be defined on $(\alpha,\beta)$. From (A6) and
\eqref{big b} we have
\begin{equation}\label{bd2}
\fb(\sigma)\delta(\sigma)^2\ \rightarrow\ \infty 
\qquad \text{as } \sigma\rightarrow
\mbox{$\alpha+$}.
\end{equation}

\begin{definition}
The partial sums of the coefficients of $\bF(s)$ and its
integral are denoted as follows:
\begin{align*}
   F(x) &\,:=\, \sum_{n \leq x}f(n) \\
   \BF(x) &\,:=\, \int_1^x F(u)du. 
\end{align*}
\end{definition}

Hayman \cite{h} 
makes a direct application of Cauchy's integral formula to express
the coefficients of a power series.
Tenenbaum \cite{te} makes a direct application of Perron's integral 
formula to express $F(x)$.
The next lemma,
where the Perron formula is used to express $\BF(x)$,
is used to derive a formula that leads to the verification of 
regular variation at infinity for $F(x)$. 
\footnote{Oppenheim  \label{footer}
(\cite{opp1}, 
\cite{opp2}) appears to state this lemma  for a particular choice
of $\bF(s)$, the zeta function connected with 
`Factorisatio Numerorum', but it is a general 
result.
On pages 210--211 of \cite{b}, one of the sources
cited for Lemma \ref{oppenheim}, the
several occurrences of $\exp(\bS(z))$ need to be replaced by
$\bS(z).$}

\begin{lemma} \label{oppenheim}
For $x \domx$ and $c>\alpha$,  
\[
\BF(x) \  = \  \frac{1}{2\pi i}\int_{c-i\infty}^{c+i\infty}
\bF(s)\;\frac{x^{s+1}}{s(s+1)}\;ds\,. 
\]
\end{lemma}
\begin{proof}
See \cite{opp1} or \Fbox{Lemma 11.22} in \cite{b}.
\end{proof}

An elementary estimate will also be needed.

\begin{lemma}\label{exp lem}
For $h,\lambda>0$ and $\kappa\in\bbR$,
\[
\int_{-h}^h e^{i\kappa-\lambda u^2}du\ =\ \sqrt{\frac{\pi}{\lambda}}e^{-\kappa^2/4\lambda}
\Big(1 + \varepsilon(h,\kappa,\lambda)\Big),
\]
where
\[
\big|\varepsilon(h,\kappa,\lambda)\big|\  <\ 
\frac{2}{h\sqrt{\lambda}}\,.
\]
\end{lemma}

The following gives the fundamental formula for $\BF(x)$.
It is this form, rather than the asymptotics that can be obtained by 
specializing 
$\sigma$ to be the saddle point $\sigma_x$, that leads to a verification of 
regular variation at infinity.

\begin{theorem}\label{main theorem}
   For $x \domx$ and $\sigma>\alpha$ 
\[
\BF(x) \ =\ \frac{x^{\sigma+1}\;\bF(\sigma)}{\sigma(\sigma+1)\sqrt{2\pi \fb(\sigma)}}
\Bigg(\exp \bigg(\frac{-\big(\fa(\sigma)+\log x\big)^2}{2\fb(\sigma)}\bigg)\, 
+\, R(x,\sigma)\Bigg)
\]
where
\[
R(x,\sigma)\ \rightarrow\ 0\qquad\text{as }\bralpha, \ \text{uniformly for }x>0.
\]
\end{theorem}

\begin{proof}
For $x \domx$ and $\sigma>\alpha$
\begin{eqnarray*}
   \BF(x)
   & = & \frac{1}{2\pi i}\int_{\sigma-i\infty}^{\sigma+i\infty}\bF(s)\;\frac{x^{s+1}}{s(s+1)}\;ds \qquad 
\text{by Lemma \ref{oppenheim}}
\\
   & = & \frac{x^{\sigma+1}}{2\pi}\int_{-\infty}^{\infty}\frac{\bF(\sigma+it)x^{it}}{(\sigma+it)(\sigma+1+it)}\;dt \\
   & = & \frac{x^{\sigma+1}}{2\pi}\big(J_1(\sigma,x) + J_2(\sigma,x) \big)
\end{eqnarray*}
where 
\begin{align*}
   J_1(\sigma,x) \ &=\  \int_{-\delta(\sigma)}^{\delta(\sigma)}\frac{\bF(\sigma+it)x^{it}}{(\sigma+it)(\sigma+1+it)}\;dt \\ 
   J_2(\sigma,x) \ &= \ \int_{|t| \geq \delta(\sigma)} 
\frac{\bF(\sigma+it)x^{it}}{(\sigma+it)(\sigma+1+it)}\;dt.\\
\intertext{Since}
\big|J_2(\sigma,x)\big| \ &\le \ \int_{|t| \geq \delta(\sigma)} 
\big|\bF(\sigma+it)\big|
\frac{dt}{\sigma^2+ t^2}
\intertext{by (A8) we immediately have}
J_2(\sigma,x) \ &=\   
\frac{\sqrt{2\pi}\;\bF(\sigma)}{\sigma(\sigma+1)\sqrt{\fb(\sigma)}}\,\ro(1)
\end{align*}
as $\bralpha$, uniformly for $x \domx$.

Let us collect some simple facts 
before estimating $J_1(\sigma,x)$.
We easily have
\begin{align} 
\frac{\sigma+1}{\sigma+1 +it}\ &=\ 1 + \ro(1)
\label{first sig}\\
\intertext{as $\bralpha$, uniformly for $|t|\le \delta(\sigma)$, since}
\frac{\sigma+1}{\sigma+1 +it}\ &=\ 1 - \frac{it}{\sigma+1+it}\nonumber\\
\intertext{and}
\Big|\frac{it}{\sigma+1+it}\Big|\ &\le\ \frac{\delta(\sigma)}{\sigma+1}\ =\ \ro(1)
\qquad\text{by \Fbox{(A4)}}.\nonumber\\
\intertext{Also}
\Big|\frac{\sigma}{\sigma +it}- 1\Big|\ &=\  \Big|\frac{it}{\sigma+it}\Big|
\ \le\ \frac{|t|}{\sigma}.
\label{second sig}
\end{align}

Let 
\[
a(\sigma,x)\ =\ \fa(\sigma)+\log x.
\]
Then by \eqref{quotient}, (A7), \eqref{first sig} and \eqref{second sig}, for 
$\sigma\in(\alpha,\beta)$ and $x>0$
\begin{eqnarray*}
J_1(\sigma,x)
&=&
\int_{|t|\le\delta(\sigma)}
\frac{\bF(\sigma+it)x^{it}} {(\sigma + it)(\sigma+it+1)} dt\\
&=&
\frac{\bF(\sigma)}{\sigma(\sigma+1)}
\int_{|t|\le\delta(\sigma)}
\exp\Big(ia(\sigma,x)t-\frac{\fb(\sigma)}{2}t^2 + \sR(\sigma+it)\Big)\\
&&\qquad\qquad\qquad\qquad\qquad\qquad\qquad \cdot\ 
\frac{\sigma(\sigma+1)}
{(\sigma + it)(\sigma+1+it)} dt\\
&=&
\frac{\bF(\sigma)}{\sigma(\sigma+1)}
\int_{|t|\le\delta(\sigma)}
\exp\Big(ia(\sigma,x)t-\frac{\fb(\sigma)}{2}t^2\Big)
\bigg(1 + \ro(1) + \rO\Big(\frac{|t|}{\sigma}\Big)\bigg)dt\\
&=&
\underbrace{
\frac{\bF(\sigma)}{\sigma(\sigma+1)}
\int_{|t|\le\delta(\sigma)}
\exp\Big(ia(\sigma,x)t-\frac{\fb(\sigma)}{2}t^2\Big)
\big(1+\ro(1)\big)
dt}_{J_{11}(\sigma,x)}\\
&& \ + \quad
\underbrace{
\frac{\bF(\sigma)}{\sigma(\sigma+1)}
\int_{|t|\le\delta(\sigma)}
\exp\Big(ia(\sigma,x)t-\frac{\fb(\sigma)}{2}t^2\Big)
\rO\Big(\frac{|t|}{\sigma}\Big)dt}_{J_{12}(\sigma,x)}.
\end{eqnarray*}

For $J_{12}(\sigma,x)$ we have, for 
$\sigma\in(\alpha,\beta)$ and $x>0$,
\begin{align*}
\big|J_{12}(\sigma,x)\big|
\ &=\ 
\rO\bigg(
\frac{\bF(\sigma)}{\sigma(\sigma+1)}
\int_{|t|\le\delta(\sigma)}
\exp\Big(-\frac{\fb(\sigma)}{2}t^2\Big)
\frac{|t|}{\sigma}dt\bigg)\\
\ &=\ 
\rO\bigg(\frac{\bF(\sigma)}{\sigma^2(\sigma+1)}
\int_0^\infty
\exp\Big(-\frac{\fb(\sigma)}{2}t^2\Big) tdt\bigg)\\
\ &=\ 
\rO\bigg(\frac{\bF(\sigma)}{\sigma(\sigma+1)\sqrt{\fb(\sigma)}}
\Big(\frac{1}{\sigma \sqrt{\fb(\sigma)}}\Big)\bigg).\\
\intertext{Thus by (A5)}
J_{12}(\sigma,x)\ &=\  
\frac{\sqrt{2\pi}\;\bF(\sigma)}{\sigma(\sigma+1)\sqrt{\fb(\sigma)}}\,\ro(1)
\end{align*}
as $\bralpha$, uniformly for $x>0$.

From Lemma \ref{exp lem} we have, for 
$\sigma\in(\alpha,\beta)$ and $x>0$,
\begin{eqnarray*}
J_{11}(\sigma,x)
&=&
\frac{\bF(\sigma)}{\sigma(\sigma+1)}
\int_{|t|\le\delta(\sigma)}
\exp\Big(ia(\sigma,x)t-\frac{\fb(\sigma)}{2}t^2\Big)
\big(1+\ro(1)\big)
dt\\
&=&
\frac{\bF(\sigma)}{\sigma(\sigma+1)}
\sqrt{\frac{\pi}{\fb(\sigma)/2}}\bigg(\exp\Big(\frac{-a(\sigma,x)^2}{2\fb(\sigma)}\Big) \\
&&\qquad\qquad\qquad\qquad\qquad +\ 
\varepsilon\big(\delta(\sigma),a(\sigma,x),\fb(\sigma)/2\big)\,+\,\ro(1) \bigg)\\
&=&
\frac{\sqrt{2\pi}\;\bF(\sigma)}{\sigma(\sigma+1)\sqrt{\fb(\sigma)}}\,
\bigg(\exp\Big(\frac{-\big(\fa(\sigma)+\log(x)\big)^2}{2\fb(\sigma)}\Big)\, +\,\ro(1)\bigg)
\end{eqnarray*}
as $\bralpha$, uniformly for $x>0$
since by Lemma \ref{exp lem} and \eqref{bd2}
\[
\Big|\varepsilon\big(\delta(\sigma),a(\sigma,x),\fb(\sigma)/2\big)\Big|\ <\ 
\frac{2}{\delta(\sigma)\sqrt{\fb(\sigma)/2}}\ =\ \ro(1).
\]


Combining these results we have
\begin{eqnarray*}
J(\sigma,x)&:=& 
J_1(\sigma,x)\ +\  J_2(\sigma,x)\\
&=&
\frac{\sqrt{2\pi}\;\bF(\sigma)}{\sigma(\sigma+1)\sqrt{\fb(\sigma)}}\,
\bigg(\exp\Big(\frac{-\big(\fa(\sigma)+\log(x)\big)^2}{2\fb(\sigma)}\Big)\, +\,\ro(1)\bigg)
\end{eqnarray*}
as $\bralpha$, uniformly for $x>0$, and the proof of the theorem is completed by
observing that
\[
\widehat{F}(x)\ =\ \frac{x^{\sigma+1}}{2\pi}\,J(\sigma,x).
\]
\end{proof}

\begin{corollary}\label{cor1}
   $\fa(\sigma)$ is strictly increasing on $(\alpha, \beta)$ 
   and as $\bralpha$
\begin{thlist}
\item
$\fa(\sigma)\  \rightarrow\  -\infty$
\item
$\displaystyle\frac{\fa(\sigma)^2}{\fb(\sigma)} \ \rightarrow\  \infty$
\item
$\big(\sigma-\alpha\big)\cdot\fa(\sigma) \rightarrow -\infty$. 
\end{thlist}
\end{corollary}

\begin{proof}
   We know that $\fa'(\sigma) = \fb(\sigma)>0$ on $(\alpha,\beta)$, so
   $\fa(\sigma)$ is strictly increasing 
   on $(\alpha, \beta)$.  Now $\BF(1) = 0$, so by Theorem 
   \ref{main theorem} with $x = 1$ we have for $\sigma\in(\alpha,\beta)$
\[
0\ =\    \frac{\bF(\sigma)}{\sigma(\sigma+1)\sqrt{2\pi \fb(\sigma)}}
   \bigg(\exp\Big(\frac{-\fa(\sigma)^2}{2\fb(\sigma)}\Big) \,+\,R(1,\sigma)\bigg)
\]
so
\[
0\ =\    \exp\Big(\frac{-\fa(\sigma)^2}{2\fb(\sigma)}\Big) \,+\,R(1,\sigma).
\]
Therefore 
$$\exp\Big(\frac{-\fa(\sigma)^2}{2\fb(\sigma)}\Big)\ \rightarrow \ 0 \baralpha,$$
and thus (b) holds.
Since $\fa(\sigma)$ decreases on $(\alpha,\beta)$ as $\bralpha$ and $\fb(\sigma)\rightarrow \infty$
by \eqref{big b} we see that (a) follows from (b).

To prove (c) we first use (a) to choose a $\gamma\in (\alpha,\beta)$ so that
$\fa(\sigma)<0$ for $\sigma\in(\alpha,\gamma)$.
Then for $\alpha <\sigma_1 < \sigma <\gamma$ we have, by the mean value theorem,
\begin{eqnarray*}
        \frac{1}{\fa(\sigma)} - \frac{1}{\fa(\sigma_1)}
        & = & -\frac{\fa'(\xi)}{\fa(\xi)^2}(\sigma-\sigma_1) \quad \text{for some $\xi \in (\sigma_1, \sigma)$} \\
        & = & -\frac{\fb(\xi)}{\fa(\xi)^2}(\sigma-\sigma_1) \\
        & = & \ro(\sigma-\sigma_1) \quad \text{as } \bralpha\\
        & = & \ro(\sigma-\alpha) \quad \text{as } \bralpha.
\end{eqnarray*}
Letting $\sigma_1 \rightarrow\mbox{$\alpha+$}$
\[
        \frac{1}{\fa(\sigma)} \ =\  \ro(\sigma-\alpha) \baralpha
\]
so $(\sigma-\alpha)\cdot\fa(\sigma) \rightarrow -\infty$ as $\bralpha$.

\end{proof}

In view of Corollary \ref{cor1}(a), from now on we will assume that $\beta$ was chosen small enough
that $\fa(\sigma)<0$ for $\sigma\in(\alpha,\beta)$.

Notice that $\fa(\sigma) \rightarrow -\infty$ 
as $\bralpha$
implies that for $x$ sufficiently large the equation
$\fa(\sigma) + \log x = 0$ has, by the continuity of $\fa(\sigma)$, a 
solution. 
In particular since $\fa'(\sigma)$ is positive on $(\alpha,\beta)$, 
for 
\[
x\,\ge\, x_0\, :=\, \exp\big(-\fa(\beta)\big)\,+\,1
\]
one has a unique solution in $(\alpha,\beta)$.

\begin{definition}
For $x\ge x_0$ (as just described)
   let $\sigma_x$ be the unique solution for $\sigma\in(\alpha,\beta)$ to the equation 
   \begin{equation}\label{the eqn}
   \fa(\sigma) + \log x \ =\  0.
\end{equation}
\end{definition}

The function $\sigma_x$ is strictly decreasing on $[x_0,\infty)$ and
\begin{equation}\label{sigma_x}
\sigma_x \rightarrow\alpha\text{\scriptsize +} \quad\text{as}\quad x\rightarrow \infty
\end{equation}
since $\displaystyle\lim_{x \rightarrow \infty}\fa(\sigma_x) = -\lim_{x \rightarrow \infty}\log x = -\infty$.

Also note that if one puts $\sigma=\sigma_x$ (where $x\ge x_0$) in the expression for
$\BF(x)$ in Theorem \ref{main theorem} then it simplifies to
\[
\BF(x) \ =\ 
\frac{x^{\sigma_x+1}\;\bF(\sigma_x)}{\sigma_x(\sigma_x+1)\sqrt{2\pi \fb(\sigma_x)}}
\big( 1 + R(x,\sigma_x)\big),
\]
where $R(x,\sigma_x)\rightarrow 0$ as $x\rightarrow \infty$. 
So we have the following.
\begin{corollary}\label{asym}
\[
 \BF(x) \ \sim\  \frac{x^{\sigma_x+1}\;\bF(\sigma_x)}{\sigma_x(\sigma_x
 +1)\sqrt{2\pi \fb(\sigma_x)}} \quad \text{as $x \rightarrow \infty$.}
\]
\end{corollary}

The choice of $\sigma=\sigma_x$ is what is commonly meant by `finding the saddlepoint',
and the resulting formula for $\BF(x)$ is the result of `applying the 
saddlepoint method'. In reality the value $s=\sigma_x$ is usually only 
near a saddle point of the integrand of the integral in Lemma \ref{oppenheim},
that is, a point $s$ where the derivative of the integrand vanishes. By
choosing the line of integration of this integral to pass through (a point near)
the saddle point one hopes to concentrate the value of the integral in a small
neighborhood of the real axis. Indeed, that is what happens for 
admissible functions. In the proof of Theorem \ref{main theorem}, the value of $\BF(x)$ 
is concentrated in the integral $J_1(x)$ when $\sigma=\sigma_x$ as $x\rightarrow \infty$,
leading to Corollary \ref{asym} above.

\begin{corollary}\label{f vs b}
As $\ralpha$
\begin{thlist}
\item
$\displaystyle
\frac{\bF(\sigma)}{\sigma\sqrt{\fb(\sigma)}} \rightarrow \infty\ $  and 
\item
$\bF(\sigma) \rightarrow \infty$. 
\end{thlist}
\end{corollary}
\begin{proof}
   Note that $\BF(2) > 0$ (since $f(1)>0$).
   Then for $\sigma\in(\alpha,\beta)$, by Theorem \ref{main theorem}
\begin{equation}\label{8.1}
\BF(2) \ =\  \frac{2^{\sigma+1}\;\bF(\sigma)}{\sigma(\sigma+1)\sqrt{2\pi \fb(\sigma)}}
\Bigg(\exp \bigg(\frac{-\big(\fa(\sigma)+\log 2\big)^2}{2\fb(\sigma)}\bigg)\, 
+\, R(2,\sigma)\Bigg).
\end{equation}
By Corollary \ref{cor1}(a) there is a $\gamma\in (\alpha,\beta)$ such that $\fa(\sigma)$
is negative on $(\alpha,\gamma)$, and thus nonzero. For $\sigma\in(\alpha,\gamma)$
we then have
\begin{eqnarray*}
\frac{\big(\fa(\sigma) + \log 2\big)^2}{\fb(\sigma)} 
& = &\frac{\fa(\sigma)^2}{\fb(\sigma)}
\bigg(1 + \frac{2\log 2}{\fa(\sigma)} + \frac{(\log 2)^2}{\fa(\sigma)^2}\bigg) .
\end{eqnarray*}
By Corollary \ref{cor1} the right hand side of this equation goes to $\infty$ as $\bralpha$,
so
\[
\exp\bigg(\frac{-\big(\fa(\sigma) + \log 2\big)^2}{2\fb(\sigma)}\bigg) \ \rightarrow\  0
\baralpha.
\]
From Theorem \ref{main theorem} we know $R(2,\sigma) \rightarrow 0$ as
$\bralpha$; and clearly
\[
   \frac{2^{\sigma+1}}{(\sigma+1)\sqrt{2\pi}} \ \rightarrow\  \frac{2^{\alpha+1}}{(\alpha+1)\sqrt{2\pi}} \ <\  \infty
\baralpha.
\]
The left side of \eqref{8.1} is a positive constant, so it follows that part (a) of this Corollary
must hold:
$\displaystyle \frac{\bF(\sigma)}{\sigma\sqrt{\fb(\sigma)}} \rightarrow \infty$
as $\bralpha.$
Then part (a) and (A5)
give 
$\bF(\sigma) \rightarrow \infty$ as $\bralpha$,
which is part (b).
\end{proof}

The next corollary shows that as $\bralpha$
we have $\bF(\sigma)$ growing much faster than any power of $\fa(\sigma)$ or 
$\fb(\sigma)$.
This leads in turn to the fact that $\bF(\sigma)$ grows much faster than any power of $\sigma-\alpha$.
Consequently $\bF(s)$ cannot have a pole at $\alpha$.

\begin{corollary}\label{waslemma2pt2}
\begin{thlist}
\item
        For all $\varepsilon > 0$, 
\[
	\fa(\sigma)\ =\ \ro\big(\bF(\sigma)^\varepsilon\big)\qquad\text{and}\qquad
       \fb(\sigma) \ =\ \ro\big(\bF(\sigma)^\varepsilon\big)\baralpha.
\]
\item
For all $r\in\mathbb{R}$, 
\[
\big(\sigma-\alpha)^r\bF(\sigma)\ \rightarrow \ \infty\baralpha.
\]
\end{thlist}
\end{corollary}

\begin{proof}
We break the proof of (a) into two claims.
\medskip

\noindent
\emph{Claim 1}:  For all $\varepsilon > 0$ and all $\gamma \in (\alpha, \beta)$ there is a $\sigma \in (\alpha, \gamma)$ such that 
\[
        \frac{|\fa(\sigma)|}{\bF(\sigma)^\varepsilon} \,<\, 1.
\]

Assume not.  Then we can choose $\varepsilon > 0$ and $\gamma \in (\alpha, \beta)$ such that for all $\sigma \in (\alpha, \gamma)$ 
\begin{equation}\label{eq1c1}
        \frac{|\fa(\sigma)|}{\bF(\sigma)^\varepsilon}\, \geq\, 1.
\end{equation}
Then for $\alpha < \sigma_1 < \sigma < \gamma$ by the mean value theorem
\begin{eqnarray*}
        \frac{1}{\bF(\sigma)^\varepsilon} - \frac{1}{\bF(\sigma_1)^\varepsilon}
        & = & \frac{-\varepsilon \bF'(\xi)}{\bF(\xi)^{1+\varepsilon}}(\sigma-\sigma_1) \quad \text{for some $\xi \in (\sigma_1, \sigma)$} \\
        & = & -\varepsilon \frac{\fa(\xi)}{\bF(\xi)^\varepsilon}(\sigma-\sigma_1) \\
        & = & \varepsilon \frac{|\fa(\xi)|}{\bF(\xi)^\varepsilon}(\sigma-\sigma_1) \\
        & \geq & \varepsilon(\sigma-\sigma_1) \quad \text{by \eqref{eq1c1}.}
\end{eqnarray*}
Letting $\sigma_1 \rightarrow\mbox{$\alpha+$}$ gives $1/\bF(\sigma)^\varepsilon \geq \varepsilon(\sigma-\alpha)$; so
\begin{equation}\label{eq2c1}
        \bF(\sigma) \ \leq\  \Big(\frac{1}{\varepsilon(\sigma-\alpha)}\Big)^{1/\varepsilon} \quad \text{for $\sigma \in (\alpha, \gamma)$}.
\end{equation}

We can also assume that $\gamma \in (\alpha, \beta)$ is such that $(\sigma-\alpha)\big|\fa(\sigma)\big| >2/\varepsilon$ for $\sigma\in (\alpha, \gamma)$ by Corollary \ref{cor1}(c).  So for $\alpha < \sigma< \sigma_2<\gamma$
\[
        -\frac{\bF'(u)}{\bF(u)}\  = \ -\fa(u)\  >\  \frac{2}{\varepsilon(u-\alpha)} \quad \text{for $u \in [\sigma,\sigma_2]$}
\]
which implies that
\[
        -\int_\sigma^{\sigma_2}\frac{\bF'(u)}{\bF(u)}du \ >\  \int_\sigma^{\sigma_2}\frac{2}{\varepsilon(u-\alpha)}du,
\]
that is,
\[
 -\big(\log \bF(\sigma_2) - \log \bF(\sigma)\big)\ >\ \frac{2}{\varepsilon}\log\Big(\frac{\sigma_2 - \alpha}{\sigma-\alpha}\Big).
\]
From this inequality and \eqref{eq2c1}
\[
        \log \bF(\sigma_2) + \frac{2}{\varepsilon}\log\Big(\frac{\sigma_2 - \alpha}{\sigma-\alpha}\Big)
        \ <\  \log \bF(\sigma) \ 
        \leq\  \frac{1}{\varepsilon}\Big(\log \frac{1}{\varepsilon} + \log \frac{1}{\sigma-\alpha}\Big).
\]
Thus
\begin{eqnarray*}
        \frac{1}{\varepsilon}\log\frac{1}{\sigma-\alpha} 
        & > & \log \bF(\sigma_2) + \frac{2}{\varepsilon}\log (\sigma_2 -\alpha) - 
\frac{1}{\varepsilon}\log\frac{1}{\varepsilon} + \frac{2}{\varepsilon}\log \frac{1}{\sigma-\alpha} \\
        & = & C + \frac{2}{\varepsilon}\log \frac{1}{\sigma-\alpha}
\end{eqnarray*}
where $C$ is independent of $\sigma$.  Hence
\[
        1\; > \;\frac{C\varepsilon}{\log \big(1/(\sigma-\alpha)\big)}\ +\ 2 \;\rightarrow\; 2 \quad \text{as } \bralpha,
\]
which is a contradiction, proving Claim 1.
\medskip

\noindent
\emph{Claim 2}:
For all $\varepsilon>0$ there is a $\gamma\in (\alpha,\beta)$ such that
\begin{equation}\label{second claim}
        \frac{|\fa(\sigma)|}{\bF(\sigma)^\varepsilon} < 1 \quad \text{for $\sigma \in (\alpha, \gamma]$.}
\end{equation}

Let $\varepsilon>0$ be given.
From Claim 1 and Corollary \ref{cor1}(b) we know that there exists a 
$\gamma \in (\alpha, \beta)$ such that 
\begin{equation}\label{choose gamma}
        \frac{|\fa(\gamma)|}{\bF(\gamma)^\varepsilon} < 1 \quad \text{and} \quad \frac{\fb(\sigma)}{\fa(\sigma)^2} < \varepsilon \quad \text{for $\sigma \in (\alpha, \gamma]$}.
\end{equation}
We will show this $\gamma$ is such that \eqref{second claim} holds.
Otherwise there is a $\sigma \in (\alpha, \gamma)$ 
such that $|\fa(\sigma)|/\bF(\sigma)^\varepsilon \geq 1$.  
By the intermediate value theorem there must be
a $\sigma \in (\alpha, \gamma)$ 
such that $|\fa(\sigma)|/\bF(\sigma)^\varepsilon = 1$.  
Letting $\sigma_1$ be the 
largest such $\sigma$ in $(\alpha,\gamma)$ we have
\[
        \frac{|\fa(\sigma_1)|}{\bF(\sigma_1)^\varepsilon} = 1 \quad \text{and} 
\quad \frac{|\fa(\sigma)|}{\bF(\sigma)^\varepsilon} < 1 \quad 
\text{for $\sigma\in (\sigma_1, \gamma]$}.
\] 
Equivalently
\begin{equation}\label{c2}
        |\fa(\sigma_1)| - \bF(\sigma_1)^\varepsilon = 0 \quad \text{and} 
\quad |\fa(\sigma)| - \bF(\sigma)^\varepsilon < 0 \quad 
\text{for $\sigma\in (\sigma_1, \gamma]$}.
\end{equation}
As $|\fa(\sigma)| = -\fa(\sigma)$ on $(\alpha, \beta)$, 
from \eqref{c2} we have
\[
        \frac{d}{d\sigma}
\Big(-\fa(\sigma) - \bF(\sigma)^\varepsilon\Big)\Big|_{\sigma = \sigma_1} \ \leq \ 0.
\]
Hence
\begin{eqnarray*}
        0
        & \leq & \fa'(\sigma_1) + \varepsilon \bF'(\sigma_1)
\bF(\sigma_1)^{\varepsilon-1} \\
        & = & \fb(\sigma_1) + \varepsilon \fa(\sigma_1)\bF(\sigma_1)^\varepsilon \\
        & = & \fb(\sigma_1) - \varepsilon \fa(\sigma_1)^2.
\end{eqnarray*}
By \eqref{choose gamma} $\fb(\sigma_1) < \varepsilon \fa(\sigma_1)^2$. 
This is a contradiction, proving Claim 2. 

From Claim 2 we immediately have
$\fa(\sigma) = \rO\big(\bF(\sigma)^{\varepsilon/2}\big)$, and thus by Corollary \ref{f vs b}(b)
$\fa(\sigma) = \ro\big(\bF(\sigma)^\varepsilon\big)$. Then from Corollary
\ref{cor1}(b)
\[
\fb(\sigma)\ =\ \ro\big(\fa(\sigma)^2\big)\ =\ \ro\big(\bF(\sigma)^\varepsilon\big)  \baralpha.
\]
This finishes the proof of (a).

Part (b) is now a trivial consequence of part (a) and Corollary \ref{cor1}(c).

\end{proof}

\begin{remark} \label{not adm}
Corollary \ref{waslemma2pt2}(b) readily shows many Dirichlet series
satisfying (A1)--(A3) are not admissible. 
\begin{thlist}
\item 
$\zeta(s)^k$, $k=1,2,\ldots$\;, 
is not
admissible as it has a pole at its abscissa $\alpha=1$.
\item
The zeta function 
\[
\prod_{j=1}^k \Big(1-n_j^{-s}\Big)^{-m_j}
\]
of a finitely generated multiplicative number system
is not admissible as it has a pole at its abscissa $\alpha=0$.
\end{thlist}
\end{remark}

\begin{corollary}\label{growth of BF}
The function $\BF(x)$ grows much faster than $x^{\alpha+1}$,
namely
\[
\lim_{x\rightarrow\infty}\frac{\BF(x)}{x^{\alpha+1}}\ =\ \infty.
\]
\end{corollary}
\begin{proof}
This is clear from 
Corollary \ref{asym},
Corollary \ref{f vs b}(a)
and \eqref{sigma_x}.
\end{proof}

\begin{definition} \label{reg var}
For $\alpha\in\bbR$,
a real-valued function $g(x)$ that is eventually
defined on the reals and eventually positive is said
to have {\bf regular variation at infinity with index 
$\alpha$}, written simply as $g(x)\in\RV_\alpha$, 
if for any $y>0$
\begin{equation}\label{yalpha}
\lim_{x\rightarrow\infty} \frac{g(xy)}{g(x)}\ =\ y^\alpha.
\end{equation}
\end{definition}

\begin{corollary} \label{BFRV}
$\BF(x) \in \RV_{\alpha+1}$.
\end{corollary}

\begin{proof}
We assume $x,y>0$.
In the expressions for $\BF(xy)$ and $\BF(x)$ given by
Theorem \ref{main theorem} let $\sigma=\sigma_x$ 
(for $x$ sufficiently large)
and divide to obtain
\begin{eqnarray}   
\frac{\BF(xy)}{\BF(x)} 
& =&  (y^{\sigma_x+1})\frac{
\exp\Big(-(\log y)^2\big/\big(2\fb(\sigma_x)\big)\Big) \,+\, R(xy,\sigma_x)}
{1+R(x,\sigma_x)}\label{xy/x}\\ 
& \rightarrow&  y^{\alpha+1} \quad\text{as }x\rightarrow \infty \nonumber
\end{eqnarray}
since both $R(x,\sigma_x) \rightarrow 0$ and 
$R(xy,\sigma_x) \rightarrow 0$ as $x\rightarrow \infty$ 
by Theorem \ref{main theorem} and \eqref{sigma_x}; and since 
$\fb(\sigma_x)\rightarrow \infty$ as 
$x\rightarrow \infty$ by \eqref{big b} and \eqref{sigma_x}.
\end{proof}

\begin{lemma} \label{drop lemma}
Let $\bG(s)\ =\ \sum _{n\ge 1} g(n)/n^s$ be a Dirichlet series
with nonnegative real coefficients
and abscissa $\alpha\ge 0$,
and let 
$\displaystyle
G(x) = \sum_{1\le n\le x} g(n)$, 
$\displaystyle\BG(x) = \int_1^x G(u)du$.

\noindent
If $\BG(x)\in\RV_{\alpha+1}$ then 
\smallskip
\begin{thlist}
\item
$G(x)\in\RV_\alpha$, and
\item
$\displaystyle
G(x)\ \sim\ \frac{\alpha+1}{x}\;\BG(x).
$
\end{thlist}
\end{lemma}
\begin{proof}
This is an immediate consequence of \Fbox{Lemma 11.21} from \cite{b}.
\end{proof}

\begin{corollary}\label{RVF}
   $F(x) \in \RV_{\alpha}$ and 
\[
  F(x) \ \sim\   
\frac{\alpha+1}{x}\;\BF(x)
 \ \sim\  \frac{x^{\sigma_x}\;\bF(\sigma_x)}{\sigma_x\sqrt{2\pi \fb(\sigma_x)}}
\]
as $x \rightarrow \infty$.
\end{corollary}

\begin{proof}
By Corollary \ref{asym}, Corollary \ref{BFRV}, Lemma \ref{drop lemma} and
\eqref{sigma_x}.
\end{proof}

\begin{corollary}\label{growth of F}
The function $F(x)$ grows much faster than $x^\alpha$,
namely
\[
\lim_{x\rightarrow\infty}\frac{F(x)}{x^\alpha}\ =\ \infty.
\]
\end{corollary}

\begin{proof}
By Corollary \ref{f vs b}(a), 
Corollary \ref{RVF} and \eqref{sigma_x}.
\end{proof}

From this Corollary it is immediate that $\zeta(s)$ is not admissible
(a fact already noted in Remark \ref{not adm}).

\section{Tenenbaum's Condtions}\label{te sect}
A version of admissibility conditions for Dirichlet series due
to Tenenbaum \cite{te}, 1988, is given in the following.\footnote{Tenenbaum actually 
uses $T(\sigma)=\sigma\fb(\sigma)/\varepsilon$, which makes our (T4) unnecessary,
and he uses the saddlepoint $\sigma_x$ instead of $\sigma$ in (T1)--(T3), which
he labels as (H2)--(H4).}

\begin{definition} \label{first te def}
Suppose $\bF(s)$ satisfies conditions {\rm (A1)--(A3)} and
there is a function
$T:(\alpha,\beta)\rightarrow (0,\infty)$, for some $\beta>\alpha$,
such that as $\sigma\rightarrow \mbox{$\alpha+$}$
\begin{enumerate}
\item[\bf (T1)]
$\sigma^2\fb(\sigma)\ \rightarrow\ \infty$
\item[\bf (T2)]
$\dfrac{\fb(\sigma)^3}{\fc(\sigma)^2} \ \rightarrow\ \infty$
\item[\bf (T3)]
$\displaystyle
T(\sigma)\;
\frac{\big|\bF(\sigma+it)\big|} {\bF(\sigma)}
\ \le \ 1$,\quad for $\sigma\in(\alpha,\beta)$ and
for $\delta(\sigma)\leq|t|\leq T(\sigma)$,\\
where $\delta(\sigma) \ =\  \big|\fb(\sigma)\fc(\sigma)\big|^{-1/5}$
\item[\bf (T4)]
$\displaystyle
\frac{\sqrt{\fb(\sigma)}}{T(\sigma)} \ \rightarrow \ 0
$
\item[\bf (T5)]
$\big|\fc(\sigma+it)\big|\le \big|\fc(\sigma)\big|$ \qquad for 
$\sigma \in (\alpha,\beta)$ and $t\in\bbR$
\item[\bf (T6)]
$\displaystyle\liminf_{\sigma\rightarrow\alpha+}\big|\fc(\sigma)\big|\ >\ 0$.
\end{enumerate}
Then we
say that $\bF(s)$ is {\bf T-admissible}, as {\bf witnessed} by $T(\sigma)$.
\end{definition}

Tenenbaum uses
$T(\sigma)= \sigma \fb(\sigma)/\varepsilon(\sigma)$ where 
$\varepsilon(\sigma)\rightarrow 0$ as $\bralpha$.
This choice of $T(\sigma)$ makes condition (T4) unnecessary.
Furthermore he gives an error term that is important to his 
applications in number 
theory, especially to the function $\psi(x,y)$. 

\begin{theorem} \label{T implies}
If $\bF(s)$ is T-admissible then it is admissible.
\end{theorem}
\begin{proof}
Let $\bF(s) = \exp\big(\bH(s)\big)$ be a T-admissible Dirichlet series as 
witnessed by $T(\sigma):(\alpha,\beta)\rightarrow(0,\infty)$.
(T1) shows that
\begin{equation} \label{repeat b}
\fb(\sigma)\rightarrow\infty \baralpha,
\end{equation}
so we can
assume that $\fb(\sigma)$ is positive on $(\alpha,\beta)$. 
From \eqref{repeat b} and (T4) it is clear that 
\[
T(\sigma)\rightarrow\infty\baralpha.
\]
By \eqref{repeat b} and (T6) one has 
\[
\delta(\sigma)\rightarrow 0\baralpha,
\]
so (A4) holds.  As (T1) is (A5)
we only need to verify that (A6)--(A8) hold. 

For (A7) we have 
for $\sigma\in(\alpha,\beta)$
and $|t|\le\delta(\sigma)$
\begin{eqnarray*}
\big|\sR(\sigma+it)\big|
&\le& 
\big|\fc(\sigma)t^3\big|\qquad \text{by \eqref{rem term}}\\
&\le& \big|\fc(\sigma)\cdot \delta(\sigma)^3\big|\\
&=& |\fc(\sigma)|\cdot\big|\fb(\sigma)\fc(\sigma)\big|^{-3/5}\\
&=& \bigg(\frac{\fc(\sigma)^2}{\fb(\sigma)^3}\bigg)^{1/5}\\
&\rightarrow&0
\qquad \text{as } \sigma\rightarrow \mbox{$\alpha+$}\quad\text{by (T2)}.
\end{eqnarray*}

For (A6) we have from (T3)
for $\sigma\in(\alpha,\beta)$
\[
T(\sigma)\;
\frac{\big|\bF\big(\sigma+i\delta(\sigma)\big)\big|}
{\bF(\sigma)}
\ \le\ 1.
\]
Multiplying this by (T4) gives
\[
\sqrt{\fb(\sigma)}\;\frac{\big|\bF\big(\sigma+i\delta(\sigma)\big)\big|}
{\bF(\sigma)} \ \rightarrow\ 0,
\]
which, in view of \eqref{abs quotient} and the fact that (A7) holds,
gives (A6).

Finally (A8) is verified as follows, where $\sigma\in(\alpha,\beta)$:
\begin{eqnarray*}
&&
\sigma\sqrt{\fb(\sigma)}\int_{|t|\ge \delta(\sigma)} \frac{\big|\bF(\sigma+it)\big|}{\bF(\sigma)}\;
\frac{dt}{\sigma^2+t^2}\\
&\le&\sigma\sqrt{\fb(\sigma)}\int_{\delta(\sigma)\le|t|\le T(\sigma)} \frac{\big|\bF(\sigma+it)\big|}{\bF(\sigma)}\; 
\frac{dt}{\sigma^2+t^2}
\quad +\quad 
2\sigma\sqrt{\fb(\sigma)}\int_{T(\sigma)}^\infty \frac{dt}{t^2}\\
&\le&\frac{\sqrt{\fb(\sigma)}}{T(\sigma)}
\bigg(\sigma\int_{\delta(\sigma)\le|t|\le T(\sigma)} 
\frac{dt}{\sigma^2+t^2}\bigg)
\ +\ 
\frac{2\sigma\sqrt{\fb(\sigma)}}{T(\sigma)}\qquad\text{by (T3)}\\
&=&
 \ro(1) \qquad \text{by (T4)}.
\end{eqnarray*}
\end{proof}

The conditions of Tenenbaum have proved to be very practical, 
giving the asymptotics for many naturally occurring examples of Dirichlet 
series to which the saddlepoint method applies. 

\begin{example} The function
\[
\bF(s)\ :=\ e^{\zeta(s)}
\]
is readily proved to be T-admissible, witnessed by $T(\sigma) = \fb(\sigma)$, 
after noting
\begin{itemize}
\item 
$\displaystyle \zeta(s)\ =\ \ \frac{1}{s-1}\,+\,g(s)$,
where $g(s)$ is holomorphic
\item
there is a constant $C>0$ such that for $\sigma\in[1,2]$ and
$|t|\ge 1$ we have
$\displaystyle \big|\zeta(\sigma+it)\big|\ \le\ C\log|t|$.
\end{itemize}

From the T-admissibility of $\exp\big(\zeta(s)\big)$ one easily has
the T-admissibility of 
\[
\bF_\lambda(s)\ :=\ \exp\big(\zeta(s-\lambda)\big)\ 
=\ \exp\Big(\sum_{n=1}^\infty n^\lambda\cdot n^{-s}\big)
\]
for $\lambda \ge 0$.\footnote{The asymptotics for $F_\lambda(x)$
are also analyzed in $\S$11.5 of \cite{b} 
by the saddlepoint method, after changing the path of the Perron
integral.
(See Footnote \ref{footer} for errata to $\S$11.5.)

R.~Warlimont \cite{warl} first pointed out the example of $\exp(\zeta(s))$ to
us. Later he found related examples of admissible
functions, such as $\exp\big(\zeta(s)^k\big)$, that subsequently turned out 
to be T-admissible as well.
}
\end{example}

\begin{example}
Tenenbaum studies the counting functions 
$\psi(x,y)$ for the zeta functions 
\[
\zeta(s,y)\ :=\ \prod_{p\le y}\big(1-p^{-s}\big)^{-1}.
\]
As noted in Remark \ref{not adm},
the functions $\zeta(s,y)$ are not admissible. 
These functions satisfy all the conditions for being T-admissible except (T3), 
and for $y$ in a suitable range (depending on $x$) they satisfy 
(T3) provided $\sigma=\sigma_x$. 
This leads to asymptotics for $\psi(x,y)$ as $x$ and $y$ 
tend to infinity with $y$ suitably constrained.
\end{example}

\begin{example}
The function
\[
\bF_k(s)\ :=\ \exp\Big(\frac{1}{1-k^{-s}}\Big)
\]
is admissible for $k=2,\ldots$\;, but not T-admissible. 
It is clear that each $\bF_k(s)$ satisfies (A1)--(A3), and has a Dirichlet 
series expansion with abscissa of convergence $\alpha=0$.
To see that $\bF_k(s)$
is not T-admissible note that
\[
\frac{\big|\bF_k(\sigma+it)\big|}{\bF_k(\sigma)}
\]
is, for each $\sigma>0$, positive and periodic as a function of $t$,
and thus does not uniformly go to 0 on $[\delta(\sigma),T(\sigma)]$ as 
$\bralpha$. Consequently 
$\bF_k(s)$ does not satisfy condition (T3).

To show that $\bF_k(s)$ is admissible let
\begin{eqnarray*}
\delta_k(\sigma)&=& (k^\sigma-1)^{7/5}\\
T_k(\sigma)&=& (k^\sigma-1)^{-3}.
\end{eqnarray*}
One has
\begin{eqnarray*}
\fa_k(\sigma)& =& -\frac{k^\sigma}{\big(k^\sigma-1\big)^2}\log k\\
\fb_k(\sigma)& =& \frac{k^{2\sigma}+k^\sigma}{\big(k^\sigma-1\big)^3}
\big(\log k\big)^2\\
\fc_k(\sigma)& =& -\frac{ k^{3\sigma}+4k^{2\sigma}+ k^\sigma }
{\big(k^\sigma-1\big)^4}
\big(\log k\big)^3.
\end{eqnarray*}

Verifying (A4)--(A6) is routine. For (A7) we proceed
as in the proof of Theorem \ref{T implies}, namely for 
$|t|\le\delta(\sigma)$ one has
\[
\big|\sR(\sigma+it)\big|\ \le\ 
\big|\fc(\sigma)\delta(\sigma)^3\big|\ \rightarrow\ 0\qquad\bralpha.
\]
This leaves (A8), which is usually the challenging part of the verification
of admissibility. First note that 
\begin{eqnarray*}
 \frac{\sigma\sqrt{\fb_k(\sigma)}}{\bF_k(\sigma)}
 \int_{|t|\ge T_k(\sigma)} \big|\bF_k(\sigma+it)\big|\;
\frac{dt}{\sigma^2+ t^2}
&\le&
 \sigma\sqrt{\fb_k(\sigma)}
 \int_{|t|\ge T_k(\sigma)} \frac{dt}{t^2}\\
&=&
 \frac{\sigma\sqrt{\fb_k(\sigma)}}{T_k(\sigma)}\ 
\rightarrow \ 0\baralpha.
\end{eqnarray*}
Thus we only need to show that
\[
 \frac{\sigma\sqrt{\fb_k(\sigma)}}{\bF_k(\sigma)}
 \int_{\delta(\sigma)\le|t|\le T_k(\sigma)} \big|\bF_k(\sigma+it)\big|\;
\frac{dt}{\sigma^2+ t^2}\ \rightarrow \ 0\baralpha.
\]
Substituting $\tau=t\log k$, we need to show 
\[
 \frac{\sigma\sqrt{\fb_k(\sigma)}}{\bF_k(\sigma)}
 \int_{\delta(\sigma)\log k\le|\tau|\le T_k(\sigma)\log k} 
 \bigg|\bF_k\Big(\sigma+i\frac{\tau}{\log k}\Big)\bigg|\;
\frac{(\log k) d\tau}{\big(\sigma\log k\big)^2+ \tau^2}\ \rightarrow \ 0
\]
as $u\rightarrow \mbox{$1+$}$.
Letting $u=k^\sigma$ it suffices to show
\[
\int_{(u-1)^{7/5}\log k}^{(u-1)^{-3}\log k}
\frac{\log u}{(u-1)^{3/2}}
\exp\bigg(\frac{2\big(\cos(\tau)-1\big)}{(u-1)\big(u^2-2u\cos(\tau)+1\big)}\bigg)
\frac{d\tau}{(\log u)^2+ \tau^2}\ \rightarrow\ 0
\]
as $u\rightarrow \mbox{$1+$}$.

One can do this by noting that as 
$u\rightarrow \mbox{$1+$}$
the integrand rapidly and uniformly approaches 0 outside neighborhoods
of radius $(u-1)^{7/5}$ about the points $\tau = 2m\pi$, indeed
much faster than $(u-1)^3$. Thus it suffices to show that
\begin{equation}\label{reduc}
\int_U
\frac{\log u}{(u-1)^{3/2}}
\exp\bigg(\frac{2\big(\cos(\tau)-1\big)}{(u-1)\big(u^2-2u\cos(\tau)+1\big)}\bigg)
\frac{d\tau}{(\log u)^2+ \tau^2}\ \rightarrow\ 0
\end{equation}
as $u\rightarrow \mbox{$1+$}$,
where $U$ is the union of the intervals 
\[
\big[2m\pi-(u-1)^{7/5},\; (2m\pi+(u-1)^{7/5}\big]
\] 
about the points $2m\pi$, $m\ge 1$, such that 
$2m\pi-(u-1)^{7/5} < (u-1)^{-3}$.
The integral in \eqref{reduc} is bounded by
\begin{equation}\label{zeta2}
2\zeta(2)
\int_{0}^{(u-1)^{7/5}}
\frac{\log u}{(u-1)^{3/2}}
\exp\bigg(\frac{2\big(\cos(\tau)-1\big)}{(u-1)\big(u^2-2u\cos(\tau)+1\big)}\bigg)
d\tau.
\end{equation}
Let $J(u,\tau)$ be the integrand in \eqref{zeta2}. Then 
\begin{eqnarray*}
\int_{0}^{(u-1)^{7/5}} J(u,\tau)d\tau
&=&
\int_{0}^{(u-1)^{3/2}} J(u,\tau)d\tau \ +\ 
\int_{(u-1)^{3/2}}^{(u-1)^{7/5}} J(u,\tau)d\tau\\
&\le& J(u,0)\cdot(u-1)^{3/2}\ +\ J\big(u,(u-1)^{3/2}\big)\cdot(u-1)^{7/5}\\
&\rightarrow& 0\qquad\text{as } u\rightarrow \mbox{$1+$}.
\end{eqnarray*}
This proves $\bF_k(s)$ is admissible, and thus the class of admissible
functions is wider than the class of T-admissible functions.

\end{example}

\section{Closure under Product} \label{closure}

The goal of this section is to prove that the product of two
admissible functions $\bF_1(s)$ and $\bF_2(s)$ with the same
abscissa of convergence is again admissible. 

\begin{theorem} \label{products}
  Suppose $\bF_1(s)$ and $\bF_2(s)$ are admissible
  with the same abscissa
  of convergence $\alpha$.
  Then $\bF_1(s)\cdot\bF_2(s)$ is admissible.
\end{theorem}

\begin{proof}
We assume 
$\bF_j(s),\delta_j(\sigma),\fb_j(s)$ 
satisfy (A1)--(A8) for $j=1,2$, and
we assume $\beta_j>\alpha$ chosen such that 
$\fb_j(\sigma)>0$ for $\sigma \in(\alpha,\beta_j)$.

Let 
\begin{eqnarray*}
\bF(s) &:=& \bF_1(s)\cdot\bF_2(s)\\
\beta&:=&\min(\beta_1,\beta_2)\\
\delta(\sigma)&:=&\min\big(\delta_1(\sigma),\delta_2(\sigma)\big)
\quad\text{for }\sigma \in (\alpha,\beta).
\end{eqnarray*}
We have 
\begin{eqnarray*}
\bF_j(s)&=&e^{\bH_j(s)}\qquad(j=1,2)\\
\bF(s)&=&e^{\bH(s)}\\
\bH(s)&=&\bH_1(s) + \bH_2(s)\\
\fa(s) &=& \fa_1(s)+\fa_2(s)\\
\fb(s) &=& \fb_1(s) + \fb_2(s)\\
\sR(s) &=& \sR_1(s) + \sR_2(s).
\end{eqnarray*}
It is easy to check that (A1)--(A5) hold for $\bF(s)$.

Next, 
\begin{eqnarray*}
  \sR(\sigma+it) 
  & = & \sR_1(\sigma+it)\,+\,\sR_2(\sigma+it) \\
  & \rightarrow & 0\qquad\text{uniformly for }|t|\le\delta(\sigma)
\end{eqnarray*}
since each of the $\sR_j$ satisfy (A7) and since $\delta(\sigma)\le\delta_j(\sigma)$
for $j=1,2$.
So (A7) also holds for $\bF$.

To prove (A6) and (A8) for $\bF$ we first observe that
for $\sigma>\alpha$, for $t\in\bbR$ and for $j=1,2$
\begin{equation} \label{peak}
\big|\bF_j(\sigma+it)\big|\,\le\, {\bF_j(\sigma)},
\end{equation}
and for $\sigma\in(\alpha,\beta)$ and $j=1,2$
\begin{eqnarray}
\fb_j(\sigma)& >& 0\\
\fb(\sigma) & =&  \fb_1(\sigma) + \fb_2(\sigma)\  \leq\  2\max\big(\fb_1(\sigma),\fb_2(\sigma)\big).
\label{max}
\end{eqnarray}
From \eqref{peak} we have 
for $\sigma>\alpha$, for $t\in\bbR$ and for $j=1,2$
\begin{equation}\label{xx}
\frac{\big|\bF(\sigma+it)\big|}{\bF(\sigma)}\ =\ 
\frac{\big|\bF_1(\sigma+it)\big|}{\bF_1(\sigma)}\cdot
\frac{\big|\bF_2(\sigma+it)\big|}{\bF_2(\sigma)}\ \le\ 
\frac{\big|\bF_j(\sigma+it)\big|}{\bF_j(\sigma)}\,.
\end{equation}

Choose $\gamma_1\in(\alpha,\beta)$ such that 
for $\sigma\in(\alpha,\gamma_1)$ and $j=1,2$
\begin{equation}\label{aa}
\fb_j(\sigma)\delta_j(\sigma)^2\ >\ 1.
\end{equation}
This is possible by \eqref{bd2}.

Now suppose that we are given $\varepsilon\in(0,1)$. 

Choose $\gamma_2\in(\alpha,\gamma_1)$ such
that for $\sigma\in(\alpha,\gamma_2)$ and $j=1,2$ 
\begin{eqnarray}
&&
\fb_j(\sigma)\cdot\exp\Big(-\fb_j(\sigma)\delta_j(\sigma)^2\Big)\ <\ \varepsilon^2
\label{bb1}\\
&&
\sigma\sqrt{\fb_j(\sigma)}\;
\int_{|t|\ge\delta_j(\sigma)}
\frac{\big|\bF_j(\sigma+it)\big|}{\bF_j(\sigma)}
\frac{dt}{\sigma^2 + t^2}
\ <\ \varepsilon.
\label{bb2}
\end{eqnarray}
We can do this because the $\bF_j$ satisfy (A6) and (A8).

Choose $\gamma\in(\alpha,\gamma_2)$ such that 
for $\sigma\in(\alpha,\gamma)$ and $j=1,2$
\begin{equation}\label{cc}
\frac{\big|\bF_j(\sigma+it)\big|}{\bF_j(\sigma)}
\ <\ 
2\exp\big(-\fb_j(\sigma)t^2/2\big)
\qquad\text{for }|t|\le\delta_j(\sigma).
\end{equation}
In view of \eqref{abs quotient} we can do this because the 
$\bF_j$ satisfy (A7).  \medskip

\noindent{\bf Claim:}
For $\sigma\in(\alpha,\gamma)$
\begin{eqnarray}
&&
\fb(\sigma)\cdot\exp\Big(-\fb(\sigma)\delta(\sigma)^2\Big)\ <\ 2\varepsilon^2
\label{goal1}\\
&&
\sigma\sqrt{\fb(\sigma)}\;
\int_{|t|\ge\delta(\sigma)}
\frac{\big|\bF(\sigma+it)\big|}{\bF(\sigma)}
\frac{dt}{\sigma^2 + t^2}
\ <\ 12\varepsilon.
\label{goal2}
\end{eqnarray}
This will prove that (A6) and (A8) hold for $\bF$.

We start by fixing $\sigma\in(\alpha,\gamma)$.

\begin{itemize}
\item[]\underline{Case (i):}\qquad
$\delta_2(\sigma)\ \le\ \delta_1(\sigma)$.\\

Then $\delta(\sigma) = \delta_2(\sigma)$.\\

\item[]\underline{Subcase (ia):} \qquad$\fb_1(\sigma)\ \le\ \fb_2(\sigma)$.\\

Then by \eqref{max} and \eqref{bb1}
\begin{eqnarray*}
\fb(\sigma)\cdot\exp\Big(-\fb(\sigma)\delta(\sigma)^2\Big)
&<& 
2\fb_2(\sigma)\cdot\exp\Big(-\fb_2(\sigma)\delta_2(\sigma)^2\Big)\\
&<&2\varepsilon^2.
\end{eqnarray*}
Also
by \eqref{max}, 
\eqref{xx} for $j=2$ and \eqref{bb2} for $j=2$ we have
\begin{eqnarray*}
&&
\sigma\sqrt{\fb(\sigma)}
\int_{|t|\ge \delta(\sigma)}
\frac{\big|\bF(\sigma+it)\big|}{\bF(\sigma)}
\frac{dt}{\sigma^2 + t^2}\\
&\le& 
\sigma\sqrt{2\fb_2(\sigma)}
\int_{|t|\ge \delta_2(\sigma_2)}
\frac{\big|\bF_2(\sigma+it)\big|}{\bF_2(\sigma)}
\frac{dt}{\sigma^2 + t^2}\ 
\\
& <& \sqrt{2}\varepsilon\ <\ 12\varepsilon.
\end{eqnarray*}

\item[]\underline{Subcase (ib):} \qquad$\fb_2(\sigma)\ \le\ \fb_1(\sigma)$.\\

By \eqref{bb1} for $j=2$
\[
\sqrt{\fb_2(\sigma)}\;
\exp\big(-\fb_2(\sigma)\delta_2(\sigma)^2/2\big)\ <\ \varepsilon,
\]
so
\begin{equation}\label{ee}
\sqrt{\fb_1(\sigma)}\;
\exp\big(-\fb_1(\sigma)\delta_2(\sigma)^2/2\big)\ <\ \varepsilon
\end{equation}
since
\[
\sqrt{x}\; \exp\big(-x\delta_2(\sigma)^2/2\big)
\]
is decreasing for $x>1/\delta_2(\sigma)^2$, 
and since by Subcase (ib) and \eqref{aa} 
for $j=2$
\[
\fb_1(\sigma)\ \ge\ \fb_2(\sigma)\ >\ \frac{1}{\delta_2(\sigma)^2}.
\]

Then by \eqref{max} and \eqref{ee}
\begin{eqnarray*}
\fb(\sigma)\cdot\exp\Big(-\fb(\sigma)\delta(\sigma)^2\Big)
&<& 
2\fb_1(\sigma)\cdot\exp\Big(-\fb_1(\sigma)\delta_2(\sigma)^2\Big)\\
&<&2\varepsilon^2.
\end{eqnarray*}

Also from \eqref{ee} we have
\[
\sqrt{\fb_1(\sigma)}\;
\exp\big(-\fb_1(\sigma)t^2/2\big)\ <\ \varepsilon
\qquad \text{for } \delta_2(\sigma)\le|t|.
\]
Combined with \eqref{cc} for $j=1$ this gives 
\begin{equation}\label{ff}
\sqrt{\fb_1(\sigma)}\; \frac{\big|\bF_1(\sigma+it)\big|}{\bF_1(\sigma)}\ <\ 2\varepsilon
\qquad \text{for } \delta_2(\sigma)\le|t|\le\delta_1(\sigma).
\end{equation}
By \eqref{max} and \eqref{xx} for $j=1$
\begin{eqnarray*}
&&
\sigma\sqrt{\fb(\sigma)}\;
\int_{|t|\ge\delta(\sigma)}
\frac{\big|\bF(\sigma+it)\big|}{\bF(\sigma)}
\frac{dt}{\sigma^2 + t^2}\\
&\le&
\sigma\sqrt{2\fb_1(\sigma)}\;
\int_{|t|\ge\delta_2(\sigma)}
\frac{\big|\bF_1(\sigma+it)\big|}{\bF_1(\sigma)}
\frac{dt}{\sigma^2 + t^2}\\
&=&
\sqrt{2}\Big(J_1(\sigma)\,+\,J_2(\sigma)\Big),
\end{eqnarray*}
where
\begin{eqnarray*}
J_1(\sigma)
&=&
\sigma\sqrt{\fb_1(\sigma)}\;
\int_{\delta_2(\sigma)\le |t|\le\delta_1(\sigma)}
\frac{\big|\bF_1(\sigma+it)\big|}{\bF_1(\sigma)}
\frac{dt}{\sigma^2 + t^2}\\
J_2(\sigma)
&=&
\sigma\sqrt{\fb_1(\sigma)}\;
\int_{|t|\ge\delta_1(\sigma)}
\frac{\big|\bF_1(\sigma+it)\big|}{\bF_1(\sigma)}
\frac{dt}{\sigma^2 + t^2}.
\end{eqnarray*}
By \eqref{ff}
\[
J_1(\sigma)\ \le\ 2\varepsilon \sigma
\int_{\delta_2(\sigma)\le |t|\le\delta_1(\sigma)}
\frac{dt}{\sigma^2 + t^2}\ <\ 2\pi\varepsilon,
\]
and by \eqref{bb2}
\[
J_2(\sigma)\ <\ \varepsilon.
\]
Thus
\[
\sigma\sqrt{\fb(\sigma)}\;
\int_{|t|\ge\delta(\sigma)}
\frac{\big|\bF(\sigma+it)\big|}{\bF(\sigma)}
\frac{dt}{\sigma^2 + t^2}\ <\ \sqrt{2}\big(2\pi+1\big)\varepsilon\ 
<\ 12\varepsilon,
\]
and the claim is proved in Case (i). 
\end{itemize}

Case (ii), where
$\delta_1(\sigma)\le\delta_2(\sigma)$, is handled likewise.
So (A6) and (A8) hold for $\bF$, and the theorem is proved.

\end{proof}

\section{Open questions}

\begin{problem} 
Is the sum of two admissible functions also admissible?
\end{problem}

\begin{problem}
Is the product of any two admissible functions also admissible?
\end{problem}

\begin{problem}
Given two admissible functions $\bF_j(x) = \exp\big(\bH_j(s)\big)$ is
the function $\exp\big(\bH_1(s)\cdot\bH_2(s)\big)$ admissible?
\end{problem}

\begin{problem}
If $\bF(s)$ is admissible, does it follow that $e^{\bF(s)}$ is also admissible?
\end{problem}

We suspect, by analogy with Hayman's work, that this is true. 

\begin{problem}
Can the notion of admissible be extended to include $\zeta(s)$?
\end{problem}

\end{document}